\documentclass[12pt,a4paper,oneside]{article}
\usepackage{authblk}
\usepackage[longnamesfirst]{natbib}
\setcitestyle{round}
\usepackage{amsmath}
\usepackage{mathrsfs}  
\usepackage{amsfonts}
\usepackage{amssymb}
\usepackage{amsthm}
\usepackage{graphicx}
\usepackage{soul, color, colortbl}
\usepackage{enumerate}
\usepackage[shortlabels]{enumitem}
\usepackage{url}
\newtheorem{assumption}{Assumption}
\newtheorem{theorem}{Theorem}

\newtheorem{lemma}{Lemma}
\newtheorem{proposition}{Proposition}
\newtheorem{remark}{Remark}
\numberwithin{assumption}{section}
\numberwithin{theorem}{section}
\numberwithin{proposition}{section}
\numberwithin{remark}{section}
\numberwithin{lemma}{section}

\newenvironment{acknowledgement}{\vspace{0.75cm} \small \noindent {\bf Acknowledgement} }%
{\small \par\noindent }

\definecolor{Gray0}{gray}{0.94}
\definecolor{Gray1}{gray}{0.82}
\definecolor{Gray2}{gray}{0.62}

\title{ {\bf Robust estimation of multivariate location and scatter in the
presence of cellwise and casewise contamination}}
\author[1]{{\bf Claudio Agostinelli}}
\author[2]{{\bf Andy Leung}}
\author[3]{{\bf Victor J. Yohai}}
\author[2]{{\bf Ruben H. Zamar}}
\affil[1]{{\small Dipartimento di Scienze Ambientali, Informatica e Statistica,
Universit\`a Ca' Foscari di Venezia,
San Giobbe, Cannaregio 873,
30121 Venezia}}
\affil[2]{{\small Department of Statistics, University of British Columbia, 3182-2207 Main Mall, Vancouver,
British Columbia V6T 1Z4, Canada}}
\affil[3]{{\small Departamento de Matem\'atica, Facultad de Ciencias Exactas y Naturales, Universidad de Buenos Aires, 
Ciudad Universitaria, Pabell\'on 1, 1426, Buenos Aires, Argentina}}

\begin{document}
\maketitle

\begin{abstract}
Multivariate location and scatter matrix estimation is a cornerstone in
multivariate data analysis. We consider this problem when the data may contain
independent cellwise and casewise outliers. Flat data sets with a
large number of variables and a relatively small number of cases are common
place in modern statistical applications. In these cases global down-weighting
of an entire case, as performed by traditional robust procedures, may lead to poor
results. We highlight the need for a new generation of robust estimators that
can efficiently deal with cellwise outliers and at the same time show 
good performance under casewise outliers.
\end{abstract}

\section{Introduction}\label{sec:INTRO}

\emph{Outliers} are a common problem for data analysts because they may have a
big detrimental effect on estimation, inference and prediction. On the other
hand, outliers could be of main interest to data analysts because they
may  represent interesting rare cases such as rocks with an unusual composition of chemical compounds and exceptional athletes in a major league. The main goal
in this article is robust estimation of multivariate location  and scatter matrix  in the presence of outliers. The estimation of these parameters is a corner stone in many applications such as principal component analysis, factor analysis, and multiple linear regression. \citet{alqallaf:2009} introduced a new contamination model where traditional
robust and affine equivariant estimators fail. To handle this new type of
outliers, we propose a new method that involves two steps: a first step of
outliers filtering, i.e., detection and replacement by missing values denoted by NA's, and a
second step of robust estimation.

\subsection*{Classical contamination model}

To fix ideas, suppose that a multivariate data set is organized in a table
with rows as cases and columns as variables, that is,  $\mathbb X = (\mathbf X_1, ..., \mathbf X_n)'$, with  $\mathbf X_i = (X_{i1},...,X_{ip})$. The vast majority of procedures for robust analysis of multivariate data are based on the classical Tukey-Huber contamination model (THCM), where a small fraction of rows in the data table may be contaminated. In THCM the contamination
mechanism is modeled as a mixture of two distributions: one corresponding to
the nominal model and the other corresponding to the outliers. More precisely,
THCM considers the following family of distributions:
\begin{equation}\label{eq:THCM}
\mathscr{H}_{\epsilon}= \{ H = (1 - \epsilon) H_0 + \epsilon \widetilde H: \widetilde H \text{ is any distribution on } \mathbb R^p\}
\end{equation}
where $H_0$ is a central parametric distribution such as the multivariate normal 
$N_p(\boldsymbol{\mu}, \mathbf{\Sigma})$ and $\widetilde H$ is an unspecified
outlier generating distribution. We then assume a case follows a distribution from the above family, that is $\mathbf X_i \sim H$ where $H \in \mathscr{H}_{\epsilon}$.  The key feature of this model is that when $\epsilon$ is small we have $\mathbf{X}_i \sim H_0$ most of the time, therefore detection and down-weighting of outlying cases makes sense and works well in practice. High breakdown point affine equivariant estimators such as MVE \citep{rousseeuw:1985}, MCD \citep{rousseeuw:1985}, S \citep{davies:1987}, MM \citep{tatsuoka:2000} and Stahel-Donoho estimators \citep{stahel:1981, donoho:1982} proceed in this general way.

\subsection*{Independent contamination model}

In many applications, however,  the contamination mechanism may be different in that
individual components (or cells) in $\mathbb X$ are independently
contaminated. This is particularly so in  the case of high dimensional data where variables
are often measured  separately  and/or obtained from different sources. For instance,
pathology and treatment information of a patient can be obtained from the
cancer registry while epidemiological information on the patients are normally
obtained through a survey. The cellwise contamination mechanism may in
principle seem rather harmless, but in fact it has far reaching consequences
including the possible breakdown of classical high breakdown point estimators.

The new contamination framework, called independent contamination model (ICM),
was presented and formalized in \citet{alqallaf:2009}. In the ICM
framework we consider a different family of distribution: 
\begin{equation}\label{eq:ICM}
\mathscr I_\epsilon = \{H: H \text{ is the distribution of } \mathbf X = (\mathbf I - \mathbf B_\epsilon) \mathbf X_0 + \mathbf B_\epsilon \widetilde{\mathbf X}  \},
\end{equation}
where $\mathbf X_0 \sim H_0$, $\widetilde{\mathbf X} \sim \widetilde H$, and $\mathbf B_\epsilon = \text{diag}(B_1,...,B_p)$, where  the $B_j$ are independent $Bin(1, \epsilon)$. 
In other words, each component of $\mathbf X$ has a probability $\epsilon$
of being independently contaminated. Furthermore, the probability $\overline{\epsilon}$ that at least one
component of $\mathbf{X}$ is contaminated is now
\[
\overline{\epsilon}=1-(1-\epsilon)^{p}.
\]
This implies that even if $\epsilon$ is small, $\overline{\epsilon}$ could be
large for large $p$, and could exceed the 0.5 breakdown point of highly  robust
affine equivariant estimators under THCM. For example, if $\epsilon=0.1$ and
$p=10$, then $\overline{\epsilon}=0.65$; if $\epsilon=0.05$ and $p=20$, then
$\overline{\epsilon}= 0.64$ and if $\epsilon=0.01$ and $p=100$, then
$\overline{\epsilon}= 0.63$.

\citet[]{alqallaf:2009} showed that for this type of contamination the
breakdown point of all the traditional $0.5$ breakdown point  and affine equivariant location
estimators  is $1 -0.5^{1/p}\to 0$ as $p \to \infty$.
It can be shown that the same holds for  robust and affine equivariant scatter estimators.
Hence we have a new manifestation  of the \emph{curse of dimensionality}: when $p$ is large, 
traditional robust estimators break down for a rather small fraction of independent contamination.

To remedy this problem, some
researchers have proposed to Winsorize potential outliers for each variable
separately. For instance, \citet{alqallaf:2002} revisited Huberized Pairwise
Covariance \citep{huber:1981}, which is constructed by using transformed
correlation coefficients calculated separately on Huberized  data  as basic building blocks. {\em Huberization} is a form of Winsorization. 
Although  pairwise robust estimators  show some robustness under ICM, they  cannot deal with THCM outliers and finely shaped  multivariate data.
Another approach to deal with ICM outliers was proposed in \citet{vanAelst:2012}. They modified the Stahel-Donoho (SD) estimator \citep{stahel:1981, donoho:1982} by calculating the  SD-outlyingness  measure and weights  on Huberized data instead of the raw data. In our simulation study   this estimator performs very well under THCM but is not sufficiently robust under ICM.
 
An alternative approach, called {\em snipping} in a recent paper by  \citet{farcomeni:2013}, consists of replacing cellwise outliers by NA.  An interesting idea introduced in \citet{farcomeni:2013} is the notion of optimizing over the snipping set. The use of snipping to fend against cellwise contamination has also been 
 suggested by other authors \citep[e.g.,][]{danilov:2010, vanAelst:2012}. \citet{farcomeni:2013} gives a procedure for clustering multivariate data where each cluster has an unknown location and scatter matrix.  This framework  can be adapted to our setting by fixing the number of clusters  to one. 
\citet{farcomeni:2013} suggested to first fix the proportion of cells in the data table to be
snipped and  then to use a maximum likelihood based procedure  to obtain  an optimal  set of snipped cells (of the same size) together with an estimate of the location and scatter matrix for each cluster. In our simulation study  this estimator performs very well under ICM but is not sufficiently robust under THCM.

A new generation of global--robust estimators that can simultaneously deal with cellwise and casewise outliers  is needed. In Section 2, we introduce a  global--robust estimator of
multivariate location and scatter. In Section 3, we show that our estimation
procedure is strongly consistent. That is, the multivariate location estimator converges a.s. to the true location and the scatter matrix estimator converges a.s. to a scalar multiple of the true scatter matrix, for a general elliptical distribution. Moreover, for a normal distribution the scalar factor is equal to one. 
In Section 4, we report the result of an extensive Monte Carlo simulation
study. In Section 5, we analyze a real data set using the proposed and several
competing estimators. In Section 6, we conclude with some remarks. Section 7 is an Appendix containing all the proofs and some additional numerical results.   

\section{Global-robust estimation under THCM and ICM}\label{sec:EST}

The main goal of this paper is to emphasize the need for robust estimation under ICM {\em and} THCM, that is, to define robust estimators that can deal with cellwise and casewise outliers.

When preprocessing multivariate data, one could try to detect cellwise outliers by applying, for instance, the ``3-sigma" rule, and replace the flagged cells by NA's. Then, an estimate of multivariate location and scatter could  be obtained using the EM-algorithm to deal with the  artificially created incomplete data.  One reason why this obvious preprocessing step is not routinely employed in multivariate robust estimation might be the lack of consistency of this procedure. Another reason might be that this approach is incapable of dealing with casewise outliers. These two limitations are addressed in our procedure by using an adaptive univariate filter \citep{gervini:2002}  followed by  Generalized S-estimator (GSE)  \citep{danilov:2012}. 

More precisely, our procedure has two steps:

\begin{enumerate}[leftmargin=*,labelindent=16pt,label=\bfseries Step \Roman*.]
\item \emph{Eliminating large cellwise outliers.} We flag cellwise
outliers and replace them
by NA's (this operation was called snipping in \citet{farcomeni:2013}).  In our case, this step  prevents cellwise contaminated cases from having large
robust Mahalanobis distances in the second step. See Section \ref{sec:EST-STEP1}.

\item \emph{Dealing with high-dimensional casewise outliers.} We apply GSE,
which has been specifically designed to deal with incomplete multivariate data with casewise outliers, to the filtered data coming from Step I. See Section
\ref{sec:EST-STEP2}.
\end{enumerate}

Full account of these steps is provided in the remaining of this section.

\subsection{Step I: Eliminating large cellwise outliers}\label{sec:EST-STEP1}

Consider a random sample of $\mathbb X = (\mathbf X_1,...,\mathbf X_n)'$, where $\mathbf X_i$ follows a distribution from $\mathscr I_\epsilon$ in (\ref{eq:ICM}).  In addition, consider a pair of initial location and dispersion estimator, $\mathbf T_{0n} = (T_{0n,1},..., T_{0n,p})$ and $\mathbf S_{0n} = (S_{0n,1},..., S_{0n,p})$. 
A common choice for $\mathbf T_{0n}$ and $\mathbf S_{0n}$  that are also adopted   in this paper are the coordinate-wise median and median absolute deviation (mad). 
%

Instead of a fixed cutoff value, we introduce an adaptive cutoff  \citep{gervini:2002} which
 is asymptotically ``correct", meaning that for clean data
the fraction of flagged outliers tends to zero as the sample size $n$ tends to
infinity. 
We identify potential outliers on each variable separately using the following 
GY-univariate filter.

We first fix a variable $(X_{1j}, X_{2j},..., X_{nj})$ and denote the standardized version of $X_{ij}$ by $Z_{ij} = (X_{ij} - T_{0n,j})/S_{0n,j}$. Let $F_j$ be a chosen reference distribution for $Z_{ij}$. An ideal choice for a reference distribution would be $F_{0j}$, the actual distribution of $(X_{ij} - \mu_{0j})/\sigma_{0j}$. Unfortunately,  the actual distribution of $Z_{ij}$ is never known in practice. Thus, we use the standard normal, $F_j = \Phi$, as a good approximation. 

The adaptive cutoff values are defined as follows. Let $\widehat{F}^+_{n,j}$ be the empirical distribution function for absolute standardized value, that is, 
\[
	\widehat{F}^+_{n,j}(t) = \frac{1}{n}\sum_{i=1}^n I( |Z_{ij}| \le t).
\]
The proportion of flagged outliers  is defined by 
\begin{equation}\label{eq:GY-d}
\begin{aligned}
d_{n,j}  &=\sup_{t\ge \eta_j} \left\{ F_j^+
(t)-\widehat{F}^+_{n,j}(t) \right\}^{+}\\ 
&  =\max_{i>i_{0}}\left\{  F_j^+(|Z|_{(i)j})-\frac{(i-1)}{n}\right\}^{+}, 
\end{aligned}
\end{equation}
where in general $\{a \}^+$ represents  the positive part of $a$ and $F^+$ is the distribution of $|Z|$ when $Z\sim F$.  Here $|Z|_{(i)j}$ is the order statistics of $|Z_{ij}|$, $i_{0}=\max\{i:|Z|_{(i)j}< \eta_j\}$, and $\eta_j = (F_j^+)^{-1}(\alpha)$ is a large quantile of $F^+$. We use $\alpha = 0.95$  throughout this paper, but other choices could be considered. Then we flag  $\lfloor nd_{n,j} \rfloor$ observations with the largest standardized value as cellwise outliers and replace them by NA's (here $\lfloor a \rfloor$ is the largest integer less than or equal to $a$).
Finally, the resulting adaptive cutoff value for $Z_{ij}$'s is 
\begin{equation}\label{eq:GY-t}
t_{n,j} =\min \left\{ t: \widehat F_{n,j}^+(t) \ge 1 - d_{n,j} \right\}, 
\end{equation}
that is, $t_{n,j} = Z_{(i_{n,j})j}$ with $i_{n,j} = n - \lfloor nd_{n,j} \rfloor$. Equivalently, we flag the $X_{ij}$'s with $|Z_{ij}| \ge t_{n,j}$.

The following proposition states that even when the actual distribution is unknown, asymptotically, 
the filter will  not wrongly flag an outlier
provided the tail of the chosen reference distribution is heavier (or equal)
 than that of the actual distribution. 

\begin{proposition}\label{prop:GY-asymptotic} 
Consider a (univariate) variable $X$ and a pair of  location and dispersion estimator $T_{0n}$ and $S_{0n}$. Suppose that $X \sim F_{0}((x-\mu)/\sigma)$ with $F_0$ continuous. If  the reference distribution $F^+$  satisfies:
\begin{equation}\label{eq:Reference distribution}
\max_{u\ge \eta}\left\{F^+(u)-F_{0}^+(u) \right\}\leq 0,
\end{equation}
 $T_{0n}\rightarrow\mu$ and $S_{0n}\rightarrow\sigma > 0$ a.s., 
then
\[
\frac{n_{0}}{n}\rightarrow0\text{ a.s.,}%
\]
where 
\[
n_{0}= \lfloor nd_{n} \rfloor.
\]
\end{proposition}

\noindent {\bf Proof:} See  the Appendix.

\subsection{Step II: Dealing with high-dimensional casewise outliers}\label{sec:EST-STEP2}

This second step introduces robustness against casewise outliers that went undetected in Step I. 
Data that emerges from Step I has \emph{holes} (i.e., NA's) that correspond to
potentially contaminated cells. To estimate the multivariate location and
scatter matrix from that data, we  use a recently developed estimator called GSE as briefly reviewed   below.

Let $\mathbf{X}_{i}=(X_{i1},...,X_{ip})^{\prime}$, $1\leq i\leq n$ be
$p$-dimensional i.i.d. random vectors that follow a distribution in an elliptical
family $\mathcal{E}(\boldsymbol{\mu}_{0},\mathbf{\Sigma}_{0})$ with density
\begin{equation}\label{eq:GSE-ellip}
f_{\mathbf X}(\mathbf{x},\boldsymbol{\mu}_{0},\mathbf{\Sigma}_{0})=\frac{1}%
{|\mathbf{\Sigma}_0|}f_0(d(\mathbf{x},\boldsymbol{\mu}_{0},\mathbf{\Sigma}%
_0)) 
\end{equation}
where $|A|$ is the determinant of $A$, $f_{0}$ is non-increasing and strictly
decreasing at 0, and
\begin{equation}\label{eq:GSE-mahasq}
d(\mathbf{x},\mathbf{m},\mathbf{C})=(\mathbf{x}-\mathbf{m})^{\prime}%
\mathbf{C}^{-1}(\mathbf{x}-\mathbf{m}) 
\end{equation}
is the squared Mahalanobis distance. We also use the normalized  squared Mahalanobis distances 
\begin{equation}\label{eq:GSE-mahasq-normalized}
d^{\ast}(\mathbf{x},\mathbf{m},\mathbf{C})=d(\mathbf{x},\mathbf{m},\mathbf{C}^{\ast}), 
\end{equation}
where $\mathbf{C}^{\ast}=\mathbf{C}/|\mathbf{C}|^{1/p},$ so $|\mathbf{C}^{\ast}|=1$.

Related  to $\mathbb X = (\mathbf X_1,..., \mathbf X_n)'$   we form the auxiliary data table of zeros and ones $\mathbb U = (\mathbf U_1,...,\mathbf U_n)'$. For  $1 \le i \le n$,  $\mathbf{U}_{i}=(U_{i1},...,U_{ip})^{\prime}$ is a $p$-dimensional random vector of zeros and ones, with ones indicating the observed entries of $\mathbf{X}_{i}$. Let $p_{i} = p(\mathbf{U}_{i})=\sum_{j=1}^{p}U_{ij}$ be the actual dimension of the observed part of $\mathbf{X}_{i}$.
Given a $p$-dimensional vector of zeros and ones $\mathbf u$, a  $p$-dimensional vector $\mathbf{m}$ and a  $p\times p$ matrix $\mathbf{A}$,
we denote by $\mathbf{m^{(\mathbf{u})}}$ and $\mathbf{A}^{(\mathbf{u})}$  the sub-vector of  $\mathbf{m}$ and the sub-matrix of $\mathbf{A}$, respectively,  with columns and rows corresponding to the positive entries in $\mathbf{u}$.  

Let $\widehat{\mathbf{\Omega}}$ be a $p \times p$ positive definite initial estimator for $\mathbf{\Sigma_0}$. 
Given the location vector $\boldsymbol{\mu}\in\mathbb{R}^{p}$ and a $p\times p$ positive definite matrix $\mathbf{\Sigma}$, we define the generalized M-scale, $s_{GS}(\boldsymbol{\mu},\mathbf{\Sigma},\widehat{\mathbf{\Omega}}, \mathbb X, \mathbb U)$, as the solution
in $s$ to the following equation:
\begin{equation}\label{eq:GSE-scale}
\sum_{i=1}^{n}c_{p(\mathbf{U}_{i})}\rho\left(  \frac{d^{\ast}\left(
\mathbf{X}_{i}^{(\mathbf{U}_{i})},\boldsymbol{\mu}^{(\mathbf{U}_{i})},\mathbf{\Sigma}^{(\mathbf{U}_{i})}\right)  }{s \, c_{p(\mathbf{U}_{i}%
)}\,\left\vert \widehat{\boldsymbol{\Omega}}^{(\mathbf{U}_{i})}\right\vert
^{1/p(\mathbf{U}_{i})}}\right)  =b\sum_{i=1}^{n}c_{p(\mathbf{U}_{i})}
\end{equation}
where $\rho(t)$ is an even, non-decreasing in $|t|$ and bounded loss function. The tuning
constants  $c_{k}$, $1\leq k\leq p$, are chosen such that
\begin{equation}\label{eq:GSE-tuning}
E_{\Phi}\left(  \rho\left(  \dfrac{||\mathbf{X}||^{2}}{c_{k}}\right)  \right)
=b,\quad\mathbf{X}\sim N_{k}(\mathbf{0},\mathbf{I}), 
\end{equation}
to ensure consistency under the multivariate normal. We consider the Tukey's bisquare rho function, $\rho(u) =\min(1, 1 - (1 - u)^{3})$, and $b = 0.5$ throughout this paper. 

The inclusion of $\widehat{\mathbf{\Omega}}$
in (\ref{eq:GSE-scale}) is needed to re-normalize the distances $d^*$ to achieve robustness. A heuristic argument for the inclusion of $\widehat{\mathbf{\Omega}}$ is as follows.
Suppose that $\widehat{\boldsymbol{\mu}} \approx \boldsymbol\mu_0$ and $\widehat{\mathbf{\Sigma}} \approx \widehat{\mathbf{\Omega}} \approx \mathbf\Sigma_0$. Then given $\mathbf U = \mathbf u$,
\[
	\dfrac{d^{*}(\mathbf X^{(\mathbf u)}, \widehat{\boldsymbol\mu}^{(\mathbf u)}, \widehat{\mathbf\Sigma}^{(\mathbf u)})}{c_{p(\mathbf u)}\left|\widehat{\mathbf{\Omega}}^{(\mathbf u)} \right|^{1/p(\mathbf u)}} \approx  \dfrac{d^{*}(\mathbf X^{(\mathbf u)}, \boldsymbol\mu_0^{(\mathbf u)}, \mathbf\Sigma_0^{(\mathbf u)})}{c_{p(\mathbf u)}\left|\mathbf\Sigma_0^{(\mathbf u)} \right|^{1/p(\mathbf u)}} \sim \dfrac{||\mathbf Y^{(\mathbf u)} ||^2}{c_{p(\mathbf u)}}
\]
where $\mathbf Y^{(\mathbf u)}$ is a $p(\mathbf u)$ dimensional random vector with an elliptical distribution. Hence, $||\mathbf Y^{(\mathbf u)} ||^2 / c_{p(\mathbf u)}$ has M-scale of 1 for the given $\rho$ function if $\mathbf Y$ is normal, and large Mahalanobis distances can be down-weighted accordingly. Here, we use extended minimum volume ellipsoid (EMVE) for $\widehat{\mathbf{\Omega}}$ as suggested in \citet{danilov:2012}. 

Generalized S-estimator is then defined by
\begin{equation}\label{eq:GSE}
(\widehat{\boldsymbol{\mu}}_{GS}, \widehat{\mathbf{\Sigma}%
}_{GS}) = \arg\min_{\boldsymbol{\mu}, \mathbf{\Sigma}} s_{GS}(\boldsymbol{\mu
}, \mathbf{\Sigma}, \widehat{\mathbf{\Omega}}, \mathbb X, \mathbb U)
\end{equation}
subject to the constraint
\begin{equation}\label{eq:GSE-constraint}
s_{GS}(\boldsymbol{\mu}, \mathbf{\Sigma},
\mathbf{\Sigma}, \mathbb X, \mathbb U) = 1.
\end{equation}

Under mild regularity assumptions, in the case of elliptical data with
$\mathbf{U}_{i}$ independent of $\mathbf{X}_{i}$ (missing completely at random
assumption) any solution to (\ref{eq:GSE}) is a consistent estimator for
the shape of the scatter matrix. Moreover, in the case of normal data, any
solution to (\ref{eq:GSE}) satisfying (\ref{eq:GSE-constraint}) is
consistent in shape and size for the \emph{true} covariance matrix. Proofs of
these claims, as well as the formulas and the derivations of the estimating
equation for GSE, can be found in \citet{danilov:2012}.

Finally our two-step location and scatter estimator is defined by 
\begin{equation}\label{eq:2SGS}
\begin{aligned}
\mathbf T_{1n} &= \widehat{\boldsymbol\mu}_{GS}( \mathbb X, \mathbb U(\mathbf t_{n}) )\\ 
\mathbf C_{1n} &= \widehat{\mathbf\Sigma}_{GS}( \mathbb X,  \mathbb U(\mathbf t_{n}))\\ 
\end{aligned}
\end{equation}
where  $\mathbf t_{n} = (t_{n,1},...,t_{n,p})$ \  ($t_{n,j}$ is defined in (\ref{eq:GY-t})) \ and 
\[
U_{ij}(t_{n,j}) = I\left( \left|\frac{X_{ij} - T_{0n,j}}{S_{0n,j}}\right| < t_{n,j}\right).
\]

\section{Consistency of GSE on filtered data}\label{sec:CONSISTENCY}

The  missing data created  in Step I is not missing at random because the missing data indicator, $\mathbb U$,  depends on the original data $\mathbb X$  (univariate outliers are declared missing). Therefore, the consistency of our two-step estimator cannot be directly derived from \citet{danilov:2012}.
However, as shown 
 in Theorem \ref{thm:consistency} below, our procedure is consistent at the central model
provided  the fraction of missing data converges to zero.
We need the following assumptions:

\begin{assumption}\label{assum:rho} 
The function $\rho$ is (i) non-decreasing in $|t|$, (ii)
strictly increasing at 0, (iii) continuous, and (iv) $\rho(0) = 0$ and (v)
$\lim_{v \to\infty} \rho(v) = 1$ (e.g. Tukey's bisquare rho function).
\end{assumption}

\begin{assumption}\label{assum:ellip} 
The random vector $\mathbf{X}$ follows a  distribution, $H_0$, in the
elliptical family defined by (\ref{eq:GSE-ellip}).
\end{assumption}

\begin{assumption}\label{assum:unique} 
Let $H_{0}$ be the distribution of $\mathbf{X}$ and
denote $\sigma(\boldsymbol{\mu},\mathbf{\Sigma})$ the solution in $\sigma$ to
the following equation
\[
E_{H_0}\left(  \rho\left(  \frac{d(\mathbf{X},\boldsymbol{\mu},\mathbf{\Sigma}%
)}{c_{p}\sigma}\right)  \right)  =b,
\]
and consider the minimization problem,
\begin{equation}
\min_{|\mathbf{\Sigma}|=1}\sigma(\boldsymbol{\mu},\mathbf{\Sigma}).
\label{eq:assum-uniqu}%
\end{equation}
We assume that (\ref{eq:assum-uniqu}) has a unique solution,
$(\boldsymbol{\mu}_{0},\mathbf{\Sigma}_{00})$, where $\mathbf{\Sigma}_{00}$ is
positive definite$.$ We also put $\sigma_{0}=\sigma(\boldsymbol{\mu}%
_{0},\mathbf{\Sigma}_{00})$.
\end{assumption}

\begin{assumption}\label{assum:eff} 
The proportion of  fully observed entries, 
\[
q_{n}=\#\{i,1\leq i\leq n:p_i = p(\mathbf{U}_{i}(\mathbf t_n))=p\}/n,
\] 
tends to one a.s. as $n$ tends to  infinity. Recall that $\mathbf t_n$  is the vector of cutoff values and $\mathbf{U}_i(\mathbf t_n)$ is the corresponding  indicator of observed entries in $\mathbf X_i$. 
\end{assumption}

\begin{remark}
\citet{davies:1987} showed that Assumption \ref{assum:ellip} implies
Assumption \ref{assum:unique}  with $\mathbf{\Sigma}_{00} = \mathbf{\Sigma}_{0}/| \mathbf{\Sigma}_{0}|$.
\end{remark}


\begin{remark}
By Proposition \ref{prop:GY-asymptotic}, the procedure described in Step I
satisfies Assumption \ref{assum:eff}, provided that the marginal distributions for the distribution that generated the data have tails which are lighter than or equally light to those of the reference distribution. That is, they satisfy  equation (\ref{eq:Reference distribution}). 
\end{remark}

\begin{theorem}\label{thm:consistency} 
Let $\mathbf{X}_{1}, ..., \mathbf{X}_{n}$ be a random sample from $H_{0}$ and $\mathbf{U}_{1}, ..., \mathbf{U}_{n}$ be as described in Section \ref{sec:EST-STEP2}. Suppose Assumptions \ref{assum:rho}--\ref{assum:eff}  hold.  
Let $(\widehat{\boldsymbol{\mu}}_{GS}, \widehat{\mathbf{\Sigma}}_{GS})$ be the GSE defined by (\ref{eq:GSE})--(\ref{eq:2SGS}). Then
\begin{enumerate}[(i)]
\item $\widehat{\boldsymbol{\mu}}_{GS} \to\boldsymbol{\mu}_{0}$ a.s. and
\item $\widehat{\mathbf{\Sigma}}_{GS}\rightarrow\sigma_{0}\mathbf{\Sigma}
_{00}$ a.s..
\item When $\mathbf{X} \sim N(\boldsymbol{\mu}_{0},\mathbf{\Sigma}_{0})$, we have
$\sigma_{0}\mathbf{\Sigma}_{00}=\mathbf{\Sigma}_{0}$.
\end{enumerate}
\end{theorem}


\noindent {\bf Proof:} See  the Appendix.


\section{Monte Carlo results}\label{sec:MCSTUDY}

We conduct a Monte Carlo simulation study to assess the performance of the proposed scatter
estimator. We consider contaminated samples from a $N_p(\boldsymbol{\mu_0},\mathbf{\Sigma_0})$ distribution. The contamination mechanisms are described below. The sample sizes are $n=100$ for
dimension $p=10$ and  $n=200$ for dimension $p=20$. 

Since the contamination models and  the estimators  considered in our simulation study
 are location and scale equivariant, we can assume without loss of
generality that the mean, $\boldsymbol\mu_0$, is equal to $\boldsymbol 0$  and the variances in  $\text{diag}(\mathbf\Sigma_0)$ are all equal to $\mathbf 1$. That is, $\mathbf\Sigma_0$ is a correlation matrix. To
account for the lack affine equivariance  of  the proposed estimator we
consider different correlation structures. For each sample in our simulation
 we create a different random correlation matrix with condition number
fixed at $CN = 100$.  Correlation matrices with high condition number are less favorable for our proposed estimator. We use the following  procedure to obtain random correlations with a fixed condition number $CN$:

\begin{enumerate}
\item For a fixed condition number CN, we first obtain a diagonal matrix
$\mathbf{\Lambda}=\mbox{diag}(\lambda_{1},...,\lambda_{p}),$  [$\lambda
_{1}<\lambda_{2}<\cdots<\lambda_{p}$] with smallest eigenvalue $\lambda_{1}=1$
and  largest eigenvalue $\lambda_{p}=\mbox{CN}$. The remaining
eigenvalues  $\lambda_{2},...,\lambda_{p-1}$ are $p-2$ sorted independent 
random variables with a uniform distribution in the interval $\left(  1,\mbox{CN}\right)  $.

\item We first generate a random $p\times p$ matrix $\mathbf{Y}$,  which elements are independent standard normal random variables. Then we form the symmetric matrix  $\mathbf{Y}^{\prime}\mathbf{Y}=\mathbf{U}%
\mathbf{V}\mathbf{U}^{\prime}$ to obtain a random orthogonal matrix $\mathbf{U}$.

\item Using the results of 1 and 2 above, we construct the random covariance matrix by $\mathbf{\Sigma}%
_{0}=\mathbf{U\Lambda U^{\prime}}$. Notice that the condition number of $\mathbf{\Sigma}_{0}$ is equal to  the desired  $CN$.

\item Convert the covariance matrix $\mathbf{\Sigma}_{0}$ into the correlation
matrix $\mathbf{R}_{0}$ as follows:
\[
\mathbf{R}_{0} =\mathbf{D}^{-1/2}\mathbf{\Sigma}_{0}\mathbf{D}^{-1/2}
\]
where
\[
 \mathbf{D} = \mbox{diag}(\sigma_{1},...,\sigma_{p}).
\]
\item  After the conversion to correlation matrix in step 4 above,  the condition number of $\mathbf{R}_{0}$  is no longer necessarily equal to $CN$. To remedy this problem, we consider the eigenvalue diagonalization   of $\mathbf{R}_{0}$

\begin{equation}\label{eq:MC-Correlation}
\mathbf{R}_{0} = \mathbf{U_0\mathbf{\Lambda}_0 U_0^{\prime}}. 
\end{equation}
where
\[
\mathbf{\Lambda}_0=\mbox{diag}(\lambda_{1}^{R_0},...,\lambda_{p}^{R_0}), \ \ \  \ \  \ \lambda_{1}^{R_0}<\lambda_{2}^{R_0} \cdots < \lambda_{p}^{R_0}.
\]
is the diagonal matrix formed using the eigenvalues of $\mathbf R_0$.
We now re-establish the desired  condition number $CN$ by redefining 
\[ \lambda_{p}^{R_0} = \mbox{CN}\times\lambda_{1}^{R_0}
\]
and using the modified eigenvalues in (\ref{eq:MC-Correlation}).

\item Repeat 4 and 5 until the condition number of $\mathbf{R}_{0}$ is within
a tolerance level (or until we reach some maximum iterations). In our
 Monte Carlo study convergence was reached after a few iteration in all the cases.
\end{enumerate}

 Two types of outliers are considered: (i)  generated by  THCM and (ii) generated by ICM. 
When the outliers are generated using THCM, we randomly replace 5\% or 10\% of the cases in the data matrix by $k\mathbf{v}$, where $k = 1,2,..., 100$ and $\mathbf{v}$ is the eigenvector corresponding to the smallest eigenvalue of $\mathbf{\Sigma}_{0}$ with length such that $\left(\mathbf{v} -\boldsymbol\mu_0\right)^{\prime}\mathbf{\Sigma}_{0}^{-1}
\left(\mathbf{v}-\boldsymbol\mu_0\right)=1$. 
Monte Carlo experiments show that the placement of outliers in this direction, $\mathbf{v}$, is the least favorable for the proposed estimator. 
When the outliers are generated using ICM, we randomly replace 5\% or 10\% of the cells in the data matrix by the value $k$ where $k=1,2,...,100$. 
The number of replicates in our simulation study is $N=500$. 

The performance of a given scatter estimator $\widehat{\mathbf{\Sigma}}$ is
measured by the Kulback-Leibler  divergence between two Gaussian distribution with the same mean and 
covariances $\mathbf{\Sigma}$ and $\mathbf{\Sigma}_{0}$:
\[
D(\mathbf{\Sigma}, \mathbf{\Sigma}_{0}) = \mbox{trace}(\mathbf{\Sigma
}\mathbf{\Sigma}_{0}^{-1}) - \log(|\mathbf{\Sigma}\mathbf{\Sigma}_{0}^{-1}|) -
p.
\]
This divergence also appears  in the likelihood ratio test statistics for testing 
the null hypothesis that a multivariate normal distribution has covariance
matrix $\mathbf{\Sigma}= \mathbf{\Sigma}_{0}$. We call this divergence measure the likelihood ratio test distance (LRT). Then the performance of an
estimator $\widehat{\mathbf{\Sigma}}$ is summarized by
\[
\overline{D}(\widehat{\mathbf{\Sigma}}, \mathbf{\Sigma}_{0}) =\frac{1}{N} \sum_{i=1}^{N}
D(\widehat{\mathbf{\Sigma}}_{i}, \mathbf{\Sigma}_{0})
\]
where $\widehat{\mathbf{\Sigma}}_{i}$ is the estimate at the $i$-th replication.

We compare the following estimators:

\begin{enumerate}[(a)]

\item MVE-S, the estimator proposed by \citet[][Section 6.7.5]{maronna:2006}.
It is an S-estimator with bisquare $\rho$ function that uses as initial value
of the iterative algorithm, an MVE estimator. The MVE estimator is computed by
subsampling with concentration step. Once the estimator of location and
covariance corresponding to one subsample are computed, the concentration step
consist in computing the sample mean and sample covariance of the [n/2]
observations with smallest Mahalanobis distance. MVE-S is implemented in the
\verb|R| package \verb|rrcov|, function
\verb|CovSest|,  option \verb|method="bisquare"|;

\item FS, the S-estimator with bisquare $\rho$ function, computed with an
iterative algorithm similar to the Fast S-estimator for regression proposed by
\citet{salibian:2006}.  FS is implemented in the
\verb|R| package \verb|rrcov|, function
\verb|CovSest|,  option  \verb|method="sfast"|;

\item MCD, the fast Minimum Covariance Determinant proposed by
\citet{rousseeuw-vandriessen:1999}  ( see also 
\citet[][Section 6.7.5]{maronna:2006} ). MCD is implemented in the
\verb|R| package \verb|rrcov|, function
\verb|CovMcd|;

\item HSD, Stahel-Donoho estimator with Huberized outlyingness proposed by
\citet{vanAelst:2012}. We use a \verb|MATLAB| code kindly provided by the
authors. The number of subsamples used in HSD is $200\times p$;

\item SnipEM, the procedure proposed in \citet{farcomeni:2013}. We use the \verb|R| code kindly
provided by the author. This method requires an initial specification of the
position of the snipped cells in the form of a binary data table. We compared (using simulation) several   possible choices for this initial set including: (a) snipping  the largest 10\% of the absolute standardized values for each variable; (b)
snipping  the largest 15\% of the absolute standardized values for each variable; and (c) snipping  the standardized values
that are more than 1.5 times the interquartile range less  the first quartile or more 
 than 1.5 times the interquartile range plus the
third quartile, for each variable. We only report the results from
case (b) as it yields the best performances.

\item 2SGS, the  two-step procedure proposed in Section
\ref{sec:EST}. This estimator is available as the \verb|TSGS| function in the \verb|R| package \verb|GSE|. 
\end{enumerate}
The tuning parameters for the  high breakdown-point estimators MVE-S, FS, and MCD are chosen to attain 0.5 breakdown point under THCM.  We have also considered pairwise scatter estimator obtained by combining bivariate S-estimator and found that this approach did not perform well in our settings (not shown here).

%
%
%
%
%

\begin{figure}
\centering
\includegraphics[scale=0.6]{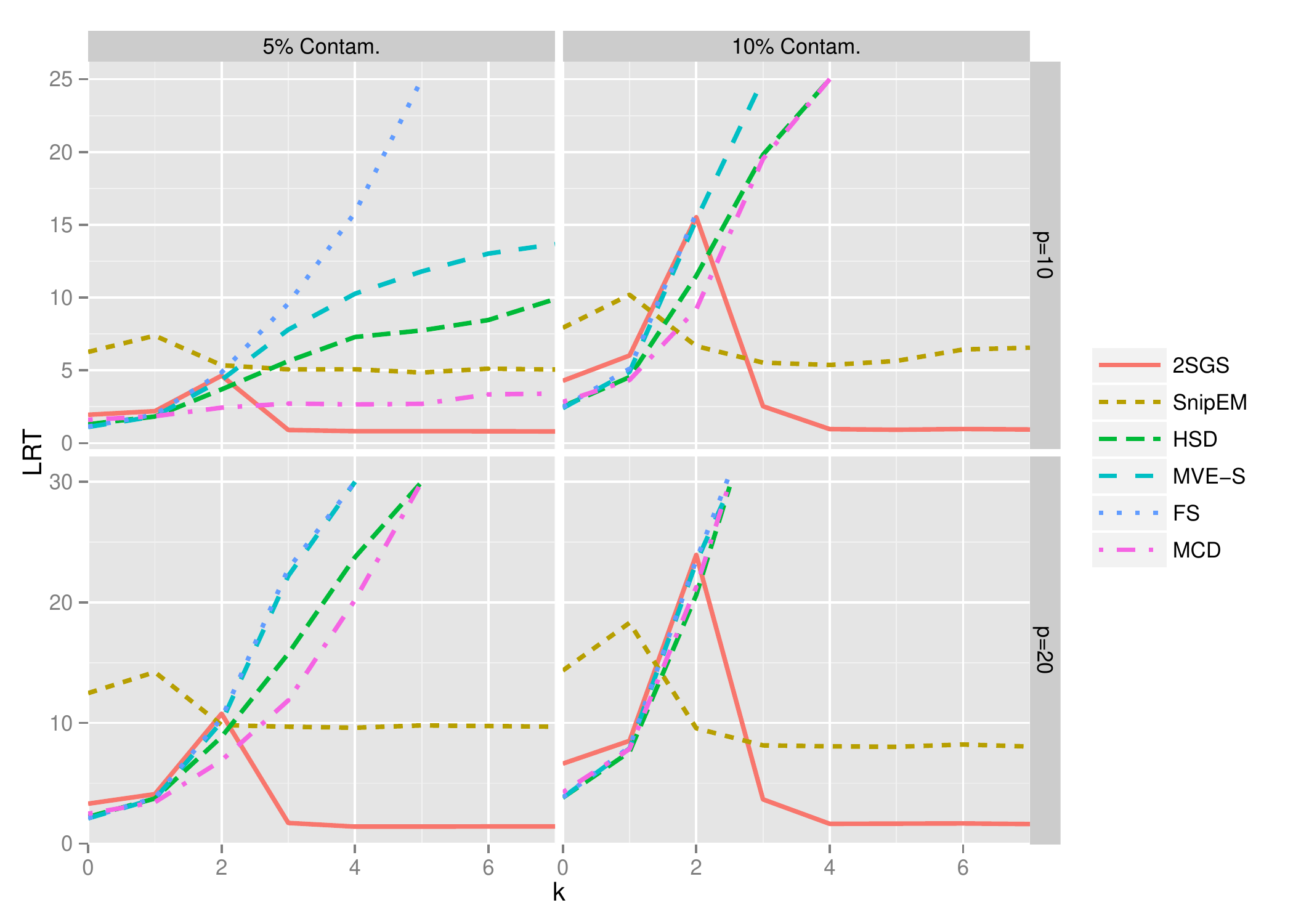}
\caption{ Average LRT distances for various contamination values, $k$, from ICM.}%
\label{fig:LRT-line-ICM}
\includegraphics[scale=0.6]{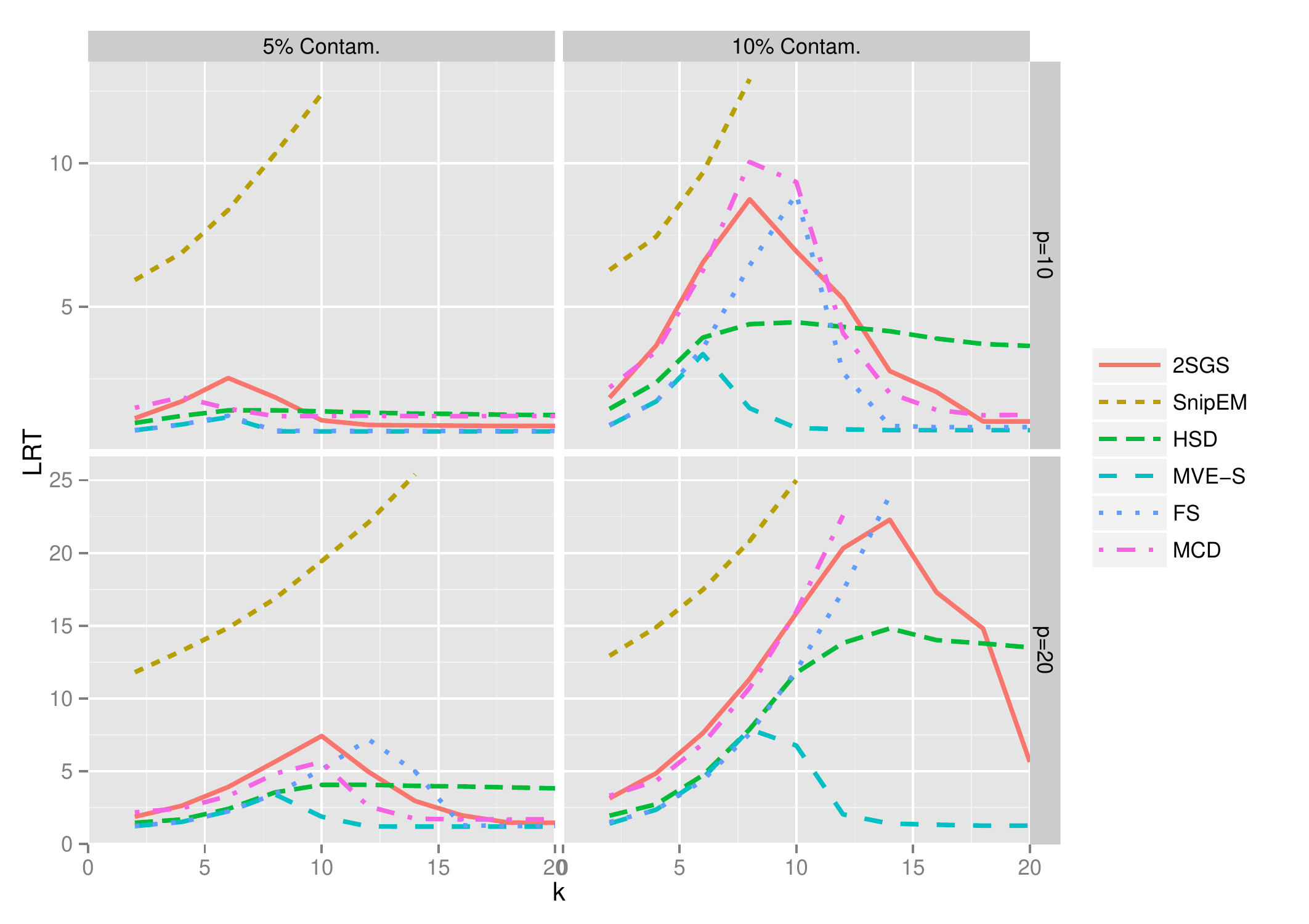}
\caption{ Average LRT distances for various contamination values, $k$, from THCM.}%
\label{fig:LRT-line-HTCM}
\end{figure}

Table \ref{tab:MCresults-robust-maxLRT} shows the maximum average LRT
distances from the true correlation matrices among the considered
contamination sizes and both contamination models. The average LRT distances behavior for different contamination sizes $k$ are displayed  in Figures \ref{fig:LRT-line-ICM} and 
\ref{fig:LRT-line-HTCM}. We notice 2SGS has the best performance under ICM. 
Not surprisingly, MVE-S has the best behavior under THCM.
However, 2SGS has an acceptable performance,  comparable with that of
main stream high breakdown point
estimators  designed for good performance under THCM.

\begin{table}[ht]
\caption{Maximum average LRT distances. Sample size is $10 \times p$. Results are
based on 500 replicates.}%
\label{tab:MCresults-robust-maxLRT}
\centering
\begin{tabular}
[c]{lrrrr|rrrr}\hline
& \multicolumn{4}{c|}{ICM} & \multicolumn{4}{c}{THCM}\\\cline{2-9}
& \multicolumn{2}{c}{Dim 10} & \multicolumn{2}{c|}{Dim 20} &
\multicolumn{2}{c}{Dim 10} & \multicolumn{2}{c}{Dim 20}\\\cline{2-9}%
Estimator & 5\% & 10\% & 5\% & 10\% & 5\% & 10\% & 5\% & 10\%\\\hline
MLE & $>$500 & $>$500 & $>$500 & $>$500 & $>$500 & $>$500 & $>$500 & $>$500 \\ 
\rowcolor{Gray1} MCD & 368.4 & $>$500 & $>$500 & $>$500 & 1.8 & 10.0 & 5.8 & 130.9 \\ 
\rowcolor{Gray1} FS & $>$500 & $>$500 & $>$500 & $>$500 & 1.2 & 8.7 & 7.2 & 204.8 \\ 
\rowcolor{Gray1} MVE-S & $>$500 & $>$500 & $>$500 & $>$500 & 1.2 & 3.3 & 3.4 & 7.9 \\ 
\rowcolor{Gray2}HSD & 11.6 & 64.7 & 75.5 & $>$500 & 1.4 & 4.5 & 4.1 & 14.8\\
\rowcolor{Gray2}SnipEM & 7.4 & 10.2 & 14.2 & 18.3 & 13.9 & 30.9 & 34.8 & 61.4 \\
\rowcolor{Gray2} 2SGS & 4.6 & 15.5 & 10.8 & 24.0 & 2.5 & 8.7 & 7.4 & 22.3 \\ \hline
\end{tabular}
\end{table}

Table \ref{tab:MCresults-efficiency-LRT} shows the finite sample relative
efficiency under clean samples for the considered robust estimates, taking the MLE average
LRT distances as the baseline. Results for larger sample sizes, not reported here, show an identical pattern, except for MCD which efficiency increases with the  sample size.

\begin{table}[htb]
\caption{Finite sample efficiency for several estimators measured by {\em relative average  LRT distances}  taking MLE as baseline. Sample size is $10 \times p$. Results are based
on 500 replicates.}%
\label{tab:MCresults-efficiency-LRT}
\centering
\begin{tabular}
[c]{lrr}\hline
Estimator & $p=10$ & $p=20$\\\hline
MLE & 1.00 & 1.00 \\ 
\rowcolor{Gray1} MCD & 0.47 & 0.66 \\ 
\rowcolor{Gray1} FS & 0.90 & 0.96 \\ 
\rowcolor{Gray1} MVE-S & 0.89 & 0.96 \\ 
\rowcolor{Gray2} HSD & 0.73 & 0.90\\
\rowcolor{Gray2} SnipEM & 0.11 & 0.28 \\
\rowcolor{Gray2} 2SGS & 0.81 & 0.84 \\ \hline
\end{tabular}
\end{table}

We also consider the barrow wheel contamination setting \citep{stahel:2009, vakili:2012} as suggested by an anonymous referee. The barrow wheel outliers are generated from a distribution that could create a large shape bias. 
The performance 2SGS is similar to the performance of the  THCM high breakdown point estimators. The results from this simulation as well as the  
computing times for our estimator (for several sample sizes and dimensions) are shown  in the Appendix

\section{Application to Chemical data}\label{sec:EXAMPLE}

We use 20 variables from a data set analyzed by \citet{smith:1984}. These
variables measure the contents (in parts per million) for 20 chemical compounds
in 53 samples of rocks in Western Australia. We compute several multivariate location and scatter estimates for this data. 

Since we suspect the occurrence of independent contamination, we
compute the $N = 53 \times 20 = 1060$ squared standardized cellwise distances and the
\[
N=53 \times 20\times 19/2=10070
\]
squared Mahalanobis distances for all the pairs $( x_{ij},x_{ik}),
i=1,2,...,53, 1\leq j<k\leq20$ using the different estimates.
 To account for multiple comparison, cellwise and pairwise distances are compared with the thresholds $(\chi
^{2}_{1})^{-1}( 0.99^{1/(np)})$ and $(\chi^{2}_{2})^{-1}( 0.99^{2/(np(p-1)})$,
respectively. To
illustrate the phenomenon of outliers propagation, full Mahalanobis distances (using all the variables) are also
computed and compared with the threshold $(\chi^{2}_{p})^{-1}( 0.99^{1/n})$.  All distances are computed using the appropriate parts from the multivariate location and scatter matrix estimates. Table \ref{tab:example-result} shows the proportion of outliers identified using the different approaches.
The proportions of identified cellwise, pairwise and casewise outliers are
 higher for robust estimators in the third generation.  In addition, the non-robust MLE flags
the smallest proportions of cellwise and pairwise outliers, and zero casewise
outliers. 

\begin{table}[htb]
\caption{Contamination summary in Chemical data based on different estimates}%
\label{tab:example-result}
\centering
\begin{tabular}
[c]{lrrr}\hline
Estimators & \multicolumn{3}{c}{Proportion of outliers}  \\\cline{2-4}
& Cell & Pair & Case  \\\hline
MLE & 0.007 & 0.008 & 0.000 \\ 
\rowcolor{Gray0}  Tyler & 0.016 & 0.024 & 0.170  \\ 
\rowcolor{Gray1}  Rocke & 0.017 & 0.027 & 0.302 \\ 
\rowcolor{Gray1}  MCD & 0.016 & 0.028 & 0.283  \\ 
\rowcolor{Gray1}   MVE & 0.024 & 0.036 & 0.283  \\ 
\rowcolor{Gray1}   FS & 0.015 & 0.027 & 0.170  \\ 
\rowcolor{Gray1}   MVE-S & 0.018 & 0.030 & 0.208  \\ 
\rowcolor{Gray2} HSDE & 0.025 & 0.038 & 0.302  \\ 
\rowcolor{Gray2} 2SGS & 0.021 & 0.033 & 0.415  \\ \hline
\end{tabular}
\end{table}

\section{Conclusions}\label{sec:CONCLUSION}

Affine equivariance,   a proven asset for achieving THCM robustness, becomes a  hindrance under ICM because of outliers propagation.

We advocate the practical and theoretical importance of ICM and point to the perils and drawbacks of relying solely
on the THCM paradigm.  ICM promotes a less aggressive cellwise down-weighting of
outliers and becomes an essential tool for modeling contamination in flat data sets 
(large in $p$ but relatively small in $n$). 
Moreover,  many 
low/moderate dimensional data sets may also be well modeled  by ICM.

We introduce a two-step procedure to  achieve  robustness under ICM and THCM. 
The first step in our procedure  is aimed at 
reducing the impact of outliers propagation and overcome the curse of
dimensionality posed by ICM. The second step  is aimed at achieving
robustness under THCM. Our procedure is not affine equivariant but nevertheless
provides fairly high resistance against both ICM and THCM outliers.
Our procedure exhibits  some  loss of robustness under THCM,  when compared with the best
performing robust affine equivariant estimators in this setting. 


We conjecture that the influence function of our estimator is the same as the influence function of the S-estimator for complete data. This conjecture is based on a similar result in \citet{gervini:2002}. They showed that the similarly derived robust regression estimator has the same influence function as the least squares estimator (they used a weighted least squares in the second step and showed that the asymptotic weights are equal to one under the central normal model). The derivation in our case seems rather involved because of the added complexity introduced by the  independent contamination model. Moreover, 
we believe that in general  the  influence function is not a very informative robustness measure. A bounded influence function is not a necessary nor sufficient condition for robustness under THCM and ICM. 

There is a need for further research on these topics.

\section{Appendix: Proofs}\label{sec:APPENDIX}
\subsection{Proof of Proposition \ref{prop:GY-asymptotic}}

Let $\widehat F_n^+$ be the empirical distribution $|Z|$ and $\widehat Z$ as defined by replacing $\mu$ and $\sigma$ with $T_{0n}$ and $S_{0n}$ respectively in the definition of $Z$.

Note that
\begin{align*}
|Z-\widehat{Z}|  &\leq  \left| \frac{X - \mu}{\sigma} - \frac{X - T_{0n}}{S_{0n}} \right| \\
&\leq  \left| \frac{X - \mu}{\sigma} - \frac{X - \mu}{S_{0n}} \right| + \frac{|T_{0n} - \mu|}{S_{0n}} \\
&  \leq \widehat{A}+\widehat{B}
\end{align*}
where $\widehat{A} \to 0$ a.s and $\widehat{B} \to 0$
a.s.. By the uniform continuity of $F^+$, given $\varepsilon>0,$ there
exists $\delta> 0$ such that $|F^+(z(1-\delta) - \delta)-F^+(z)|\leq\varepsilon/2$. With probability one there exists $n_1$ such that
$n\geq n_1$ implies $|\widehat{A}|$ $<\delta$ and $|\widehat{B}|<\delta$. By the Glivenko-Cantelli Theorem, with probability one there
exists $n_2$ such that $n\geq n_2$ implies that $\sup_{z}|\widehat F_n^+(z) -F^+(z)|\leq\varepsilon/2$. Let $n_3=\max(n_1,n_2)$, then
$n\geq n_3$ imply
\begin{align*}
\widehat{F}_{n}^+(z)  &  \geq \widehat F_{n}^+(z(1-\delta) -\delta)\\
&  =\left( \widehat F_{n}^+(z(1-\delta) -\delta)-F_0^+(z(1-\delta
)-\delta)\right) \\
&  \qquad+(F_0^+(z(1-\delta) -\delta)-F_0^{+}(z))+(F_0^{+}(z)- F^+(z))+F^{+}(z)
\end{align*}
and then
\begin{align*}
\sup_{z> \eta}( F^+(z)-\widehat{F}_{n}^+(z))  &
\leq\sup_{z> \eta}\left| F_0^+(z(1-\delta) -\delta)- \widehat F_{n}^{+}(z(1-\delta)
-\delta)\right| \\
& \qquad+\sup_{z> \eta}\left| F_0^{+}(z(1-\delta)-\delta)-F_0^{+}(z)\right| \\
& \qquad+\sup_{z> \eta}(F^{+}(z)-F_0^{+}(z))\\
&  \leq\varepsilon
\end{align*}
This implies that $n_{0}/n\rightarrow0$ a.s..

\subsection{Proof of Theorem \emph{\ref{thm:consistency}}}

We need the following Lemma proved in \citet{yohai:1985}.

\begin{lemma}
\label{lemma:consistency} Let $\{\mathbf{Z}_{i}\}$ be i.i.d. random vectors
taking values in $\mathbb{R}^{k}$, with common distribution $Q$. Let
$f:\mathbb{R}^{k}\times\mathbb{R}^{h}\rightarrow\mathbb{R}$ be a continuous
function and assume that for some $\delta>0$ we have that
\[
E_Q\left[  \sup_{||\lambda-\lambda_{0}||\leq\delta}|f(\mathbf{Z},\lambda
)|\right]  <\infty.
\]
Then if $\widehat{\lambda}_{n}\rightarrow\lambda_{0}$ a.s., we have
\[
\frac{1}{n}\sum_{1=1}^{n}f(\mathbf{Z}_{i},\widehat{\lambda}_{n})\rightarrow
E_Q\left[  f(\mathbf{Z},\lambda_{0})\right]  \text{ a.s..}%
\]
\end{lemma}

\emph{Proof of Theorem \ref{thm:consistency},}

Define
\begin{equation}\label{cons0}
(\widehat{\boldsymbol{\mu}}_{GS},\widetilde{\mathbf{\Sigma}}_{GS})=\arg
\min_{\boldsymbol{\mathbf{\mu}},|\mathbf{\Sigma}|=1}s_{GS}%
(\boldsymbol{\mathbf{\mu}},\mathbf{\Sigma},\widehat{\mathbf{\Omega}}).
\end{equation}
We drop out $\mathbb X$ and $\mathbb U$ in the argument to simplify the notation. 
Since $s_{GS}(\boldsymbol{\mu},\lambda\mathbf{\Sigma},\widehat{\mathbf{\Omega
}})=s_{GS}(\boldsymbol{\mu},\mathbf{\Sigma},\widehat{\mathbf{\Omega}})$, to prove Theorem \ref{thm:consistency} it is
enough to show
\begin{description}\item[(a)]
\begin{equation}\label{cons-1}
(\widehat{\boldsymbol{\mu}}_{Gs},\ \widetilde{\mathbf{\Sigma}}_{GS}%
)\rightarrow(\boldsymbol{\mu}_{0},\mathbf{\Sigma}_{00})\text{ a.s., \quad \quad and}
\end{equation}
\end{description}
\begin{description}\item[(b)] 
\begin{equation}\label{cons00}
s_{GS}(\widehat{\boldsymbol{\mathbf{\mu}}}_{GS},\widetilde{\mathbf{\Sigma}%
}_{GS},\widetilde{\mathbf{\Sigma}}_{GS})\rightarrow\sigma_{0}\text{ a.s.}.
\end{equation}
\end{description}

Note that since we have
\[
E_{H_0}\left(  \rho\left(  \frac{d\left(  \mathbf{X}
,\boldsymbol{\mathbf{\mu}}_{0},\mathbf{\Sigma}_{0}\right)  }{\sigma_{0}
c_{p}\,}\right)  \right)  =b,
\]
 then part (i) of Lemma 6 in the Supplemental Material of \citet{danilov:2012} implies that given $\varepsilon>0,$ there exists $\delta > 0$ such that 
\begin{equation}\label{cons1}
\underset{n\rightarrow\infty}{\underline{\lim}}\inf_{(\boldsymbol{\mu},\mathbf{\Sigma
})\in C_{\varepsilon}^{C},|\mathbf{\Sigma}|=1}\frac{1}{n}\sum_{i=1}^{n}%
c_{p}\rho\left(  \frac{d\left(  \mathbf{X}_{i},\boldsymbol{\mathbf{\mu}%
},\mathbf{\Sigma}\right)  }{\sigma_{0}c_{p}\,(1+\delta)}\right)
>(b+\delta)c_{p}, 
\end{equation}
where $C_{\varepsilon}$ is a neighborhood of $(\boldsymbol{\mu}_{0},\mathbf{\Sigma
}_{00})$ of radius $\varepsilon$ and if $A$ is a set, then $A^{C}$ denotes its
complement. In addition, by part (iii) of the same Lemma we have for any $\delta > 0$,
\begin{equation}\label{cons2}
\lim_{n\rightarrow\infty}\frac{1}{n}\sum_{i=1}^{n}c_{p}\rho\left(
\frac{d\left(  \mathbf{X}_{i},\boldsymbol{\mathbf{\mu}}_{0},\mathbf{\Sigma
}_{00}\right)  }{\sigma_{0}c_{p}\,(1+\delta)}\right)  <b\,c_{p}. 
\end{equation}

Let 
\[
Q_i(\boldsymbol\mu, \mathbf \Sigma) = c_p \rho\left(  
\frac{d\left(  \mathbf X_i,\boldsymbol{\mu},\mathbf{\Sigma}\right)}
{\sigma_0  c_p (1 + \delta)}\right)
\]
and 
\[
Q_i^{(\mathbf U)}(\boldsymbol\mu, \mathbf \Sigma) = c_{p(\mathbf U_i)} \rho\left(  
\frac{d^*\left(  \mathbf X_i^{(\mathbf U_i)},\boldsymbol{\mu}^{(\mathbf U_i)},\mathbf{\Sigma}^{(\mathbf U_i)}\right)}
{S \,c_{p(\mathbf U_i)}\, \left|\widehat{\boldsymbol\Omega}^{(\mathbf U_i)}\right| ^{1/p(\mathbf U_i)}}\right),
\]
Now if $|\mathbf{\Sigma}|=1$ and $S=\sigma_{0}(1+\delta
)/|\widehat{\mathbf{\Omega}}|^{1/p}$,
we have
\begin{equation}\label{cons3}
\begin{aligned}
\frac{1}{n}\sum_{i=1}^{n} Q_i^{(\mathbf U)}(\boldsymbol\mu, \mathbf \Sigma)   &=\frac
{1}{n}\sum_{p_{i}=p} Q_i(\boldsymbol\mu, \mathbf \Sigma)  
  +\frac{1}{n}\sum_{p_{i}\neq p} Q_i^{(\mathbf U)}(\boldsymbol\mu, \mathbf \Sigma) . 
\end{aligned}
\end{equation}

We also have%

\begin{equation}\label{cons4}
\frac{1}{n}\sum_{p_{i}\neq p}Q_i^{(\mathbf U)}(\boldsymbol\mu, \mathbf \Sigma) \leq
c_{p}(1-t_{n}) 
\end{equation}
and therefore by Assumption \ref{assum:eff} we have%
\begin{equation}\label{cons41}
\lim_{n\rightarrow\infty}\sup_{\boldsymbol{\mu},|\mathbf{\Sigma}|=1}\frac
{1}{n}\sum_{p_{i}\neq p} Q_i^{(\mathbf U)}(\boldsymbol\mu, \mathbf \Sigma) =0 \text{ a.s..} 
\end{equation}
Similarly we can prove that
\begin{equation}
\lim_{n\rightarrow\infty}\sup_{\boldsymbol{\mu},|\mathbf{\Sigma}|=1}\frac
{1}{n}\sum_{p_{i}\neq p} Q_i(\boldsymbol\mu, \mathbf \Sigma) = 0\text{
a.s.} \label{cons42}%
\end{equation}
and
\begin{equation}\label{cons45}
c_{p}-\frac{1}{n}\sum_{i=1}^{n}c_{p(\mathbf{U}_{i})}\rightarrow0,\text{ a.s..}
\end{equation}
Then, from (\ref{cons1}) and (\ref{cons3})--(\ref{cons45}) we get
\begin{equation}\label{cons5}
\underset{n\rightarrow\infty}{\underline{\lim}}\inf_{(\boldsymbol{\mu},\mathbf{\Sigma
})\in C_{\varepsilon}^{C},|\mathbf{\Sigma}|=1}\frac{1}{n}\sum_{i=1}%
^{n} Q_i^{(\mathbf U)}(\boldsymbol\mu, \mathbf \Sigma)   >  
(b+\delta)\lim_{n\rightarrow\infty}\frac{1}{n}\sum_{i=1}^{n}c_{p(\mathbf{U}_{i})}  =(b+\delta)c_{p}\ \text{a.s..} 
\end{equation}
Using similar arguments, from (\ref{cons2}) we can prove
\begin{equation}\label{cons6}
\lim_{n\rightarrow\infty}\frac{1}{n}\sum_{i=1}^{n}
Q_i^{(\mathbf U)}(\boldsymbol\mu_0, \mathbf \Sigma_{00})
  <  b\lim_{n\rightarrow\infty}\frac{1}{n}
\sum_{i=1}^{n}c_{p(\mathbf{U}_{i})} =b \, c_{p}\text{ a.s..} 
\end{equation}
Equations (\ref{cons5})--(\ref{cons6}) imply that
\[
\underset{n\rightarrow\infty}{\underline{\lim}}\inf_{(\boldsymbol{\mu},\mathbf{\Sigma
})\in C_{\varepsilon}^{C},|\mathbf{\Sigma}|=1}s_{GS}(\boldsymbol{\mathbf{\mu}%
},\mathbf{\Sigma},\widehat{\mathbf{\Omega}})>S\text{ a.s.}%
\]
and%
\[
\lim_{n\rightarrow\infty}s_{GS}(\boldsymbol{\mathbf{\mu}}_{0},\mathbf{\Sigma
}_{00},\widehat{\mathbf{\Omega}})<S\text{ a.s..}%
\]
Therefore, with probability one there exists $n_{0}$ such that for $n>n_{0}$
we have \ $(\widehat{\boldsymbol{\mu}}_{GS},\widetilde{\mathbf{\Sigma}}%
_{GS})\in$ $C_{\varepsilon}^{C}$. Then $(\widehat{\boldsymbol{\mu}}%
_{GS},\widetilde{\mathbf{\Sigma}}_{GS})\rightarrow(\boldsymbol{\mu}%
_{0},\mathbf{\Sigma}_{00})$ a.s. proving (a). \newline

Let 
\[
P_i( \boldsymbol\mu, \mathbf\Sigma, s) = c_p \rho\left(  \frac{d\left(
\mathbf{X}_{i},\boldsymbol{\mu},\mathbf{\Sigma}\right)  }{c_p \,\,s}\right)
\]
and 
\[
P_i^{(\mathbf U)}( \boldsymbol\mu, \mathbf\Sigma, s) = c_{p(\mathbf{U}_{i})} \rho\left(  \frac{d\left(
\mathbf{X}_{i}^{(\mathbf{U}_{i})},\boldsymbol{\mu}
^{(\mathbf{U}_{i})},\mathbf{\Sigma}^{(\mathbf{U}_{i})}\right)  }{c_{p(\mathbf{U}_{i})}\,\,s}\right).
\]
Since $|\widetilde{\mathbf{\Sigma}}_{GS}|=1$, we have that
$s_{GS}(\widehat{\boldsymbol{\mu}}_{GS},\widetilde{\mathbf{\Sigma}}%
_{GS},\widetilde{\mathbf{\Sigma}}_{GS})$ is the solution in $s$ in the
following equation%

\begin{equation}\label{cons7}
\frac{1}{n}\sum_{i=1}^{n} P_i^{(\mathbf U)}( \widehat{\boldsymbol\mu}_{GS}, \widetilde{\mathbf\Sigma}_{GS}, s)  =\frac{b}{n}\sum_{i=1}
^{n}c_{p(\mathbf{U}_{i})}.
\end{equation}
Then, to prove (\ref{cons00}) it is enough to show that for all
$\varepsilon>0$
\begin{equation}\label{cons8}
\begin{aligned}
&\lim_{n\rightarrow\infty}\frac{1}{n}\sum_{i=1}^{n}
P_i^{(\mathbf U)}( \widehat{\boldsymbol\mu}_{GS}, \widetilde{\mathbf\Sigma}_{GS}, \sigma_0 + \varepsilon) <b \, c_{p}\text{ a.s. \quad and} \\
&\lim_{n\rightarrow\infty}\frac{1}{n}\sum_{i=1}^{n}
P_i^{(\mathbf U)}( \widehat{\boldsymbol\mu}_{GS}, \widetilde{\mathbf\Sigma}_{GS}, \sigma_0 - \varepsilon) > b \, c_{p}\text{ a.s.}
\end{aligned}
\end{equation}
Using Assumption \ref{assum:eff}, to prove (\ref{cons8}) it
is enough to show
\begin{equation} \label{cons9}
\begin{aligned}
&\lim_{n\rightarrow\infty}\frac{1}{n}\sum_{i=1}^{n}
P_i( \widehat{\boldsymbol\mu}_{GS}, \widetilde{\mathbf\Sigma}_{GS}, \sigma_0 + \varepsilon) <b \, c_{p}\text{ a.s. \quad and} \\
&\lim_{n\rightarrow\infty}\frac{1}{n}\sum_{i=1}^{n}
P_i( \widehat{\boldsymbol\mu}_{GS}, \widetilde{\mathbf\Sigma}_{GS}, \sigma_0 - \varepsilon) > b \, c_{p}\text{ a.s.}
\end{aligned}
\end{equation}
It is immediate that%
\[
E\left(  \rho\left(  \frac{d\left(  \mathbf{X},\boldsymbol{\mu}%
_{0},\mathbf{\Sigma}_{0}\right)  }{c_{p}\,(\sigma_{0}+\varepsilon)}\right)
\right)  <E\left(  \rho\left(  \frac{d\left(  \mathbf{X},\boldsymbol{\mu
}_{0},\mathbf{\Sigma}_{0}\right)  }{c_{p}\,\sigma_{0}}\right)  \right)  =b
\]
and%
\[
E\left(  \rho\left(  \frac{d\left(  \mathbf{X},\boldsymbol{\mu}%
_{0},\mathbf{\Sigma}_{0}\right)  }{c_{p}\,(\sigma_{0}-\varepsilon)}\right)
\right)  >E\left(  \rho\left(  \frac{d\left(  \mathbf{X},\boldsymbol{\mu
}_{0},\mathbf{\Sigma}_{0}\right)  }{c_{p}\,\sigma_{0}}\right)  \right)  =b.
\]
Then equations (\ref{cons9}) follow from Lemma
\ref{lemma:consistency} and part (a). This proves (b).

\subsection{Investigation on the performance on the barrow wheel outliers}

\begin{figure}[h]
\centering
\includegraphics[scale=0.6]{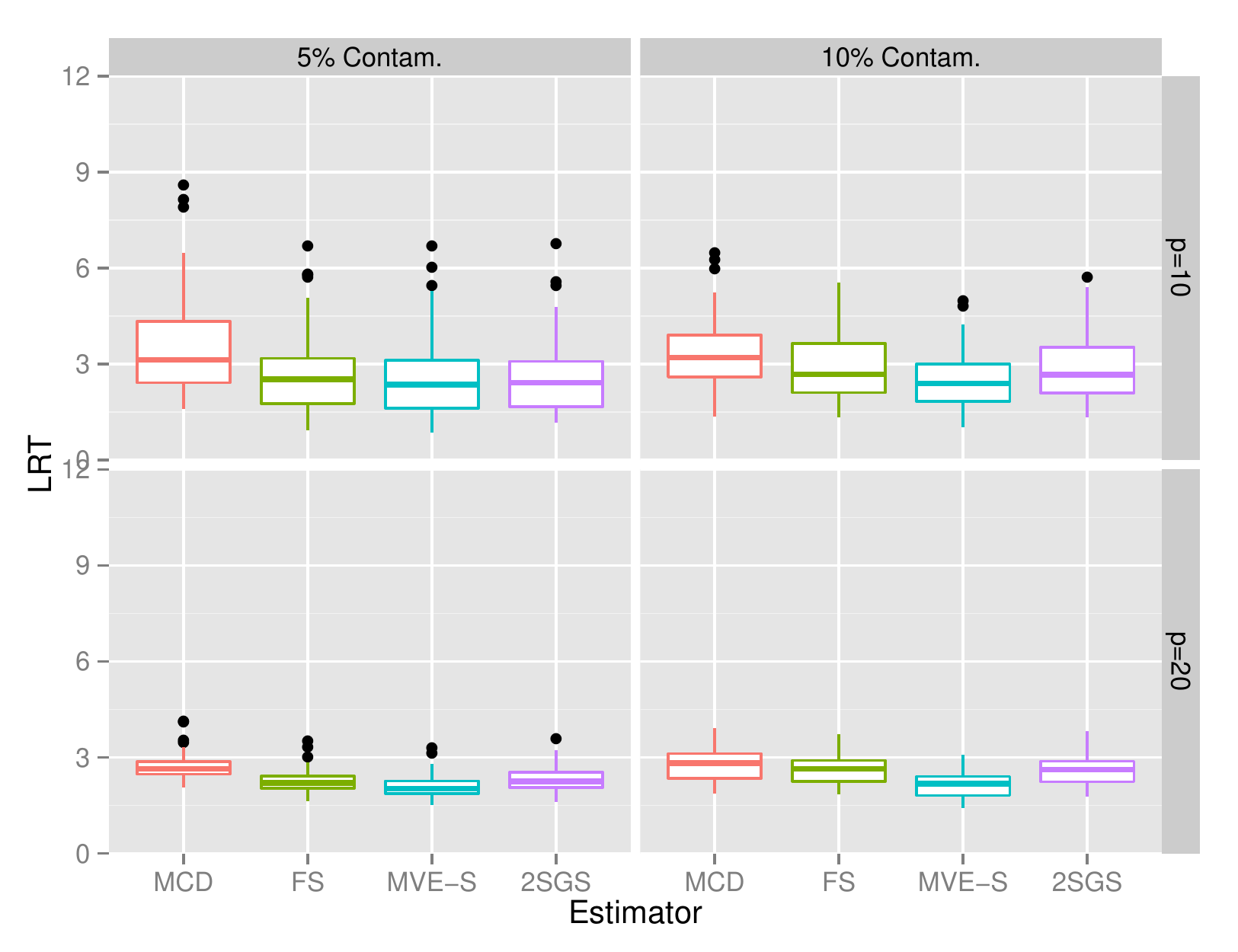}
\caption{LRT distances under barrow-wheel contamination setting. }
\label{fig:barrow-wheel}
\end{figure}

An anonymous referee suggested considering the performance of 2SGS under the barrow wheel contamination setting \citep{stahel:2009, vakili:2012}. We conduct a Monte Carlo study to compare the performance of 2SGS with three second generation estimators under 5\% and 10\% of outliers from the barrow wheel distribution. The data is generated using 
the  R package \texttt{robustX} with default parameters. The three second generation estimators are: the fast Minimum Covariance Determinant (MCD), the fast S-estimator (FS), and the S-estimator (S), described in Section \ref{sec:MCSTUDY}.  The sample size are $n = 10\times p$, for $p=10$ and $20$. The results in terms of the LRT measure are graphically displayed in Figure \ref{fig:barrow-wheel}.

\subsection{Timing experiment}

Table \ref{tab:timing} shows the mean time needed to compute 2SGS for data with  cellwise or casewise outliers as described in Section 5. We consider 10\% contamination and several sample sizes and dimensions. We use the random correlation structures as described in Section 4. For each pair of dimension and sample size, we average the computing times over 250 replications for each of the following setups: (a) cellwise contamination with $k$  generated from $U(0,6)$ and (b) casewise contamination with $k$  generated from $U(0,20)$.

\begin{table}[ht]
\caption{Average ``CPU time" -- in seconds of a 2.8 GHz Intel Xeon -- evaluated using the \texttt{R} command, \texttt{system.time}.} \label{tab:timing}
\centering
\begin{tabular}{llrr}
  \hline
$p$ & $n$ & Cellwise & Casewise \\ 
  \hline
5 & 50 & 0.03 & 0.03 \\ 
   & 100 & 0.04 & 0.04 \\ 
  10 & 100 & 0.12 & 0.10 \\ 
   & 200 & 0.17 & 0.13 \\ 
  15 & 150 & 0.40 & 0.28 \\ 
   & 300 & 0.60 & 0.40 \\ 
  20 & 200 & 1.03 & 0.73 \\ 
   & 400 & 1.88 & 1.06 \\ 
  25 & 250 & 2.52 & 1.62 \\ 
   & 500 & 4.58 & 2.45 \\ 
  30 & 300 & 5.08 & 3.26 \\ 
   & 600 & 8.47 & 5.16 \\ 
  35 & 350 & 9.30 & 6.13 \\ 
   & 700 & 15.64 & 9.79 \\ 
   \hline
\end{tabular}
\end{table}

\begin{acknowledgement}
Victor Yohai research was partially supported by Grants W276 from Universidad of
Buenos Aires,PIP 112-2008-01-00216 and
112-2011-01-00339 from CONICET and PICT2011-0397 from ANPCYT, Argentina.
Ruben Zamar and Andy Leung research were partially funded by the Natural Science and Engineering Research Council of Canada.
\end{acknowledgement}

\bibliographystyle{hplainnat}
\bibliography{FilterGS}

\end{document}